\documentclass[11pt,leqno]{article}
\usepackage{amssymb,latexsym,amsmath,theorem}
\newtheorem{theorem}{Theorem}
\newtheorem{proposition}{Proposition}
\newtheorem{corollary}{Corollary}
{\theorembodyfont{\rmfamily} \newtheorem{definition}{Definition}}
{\theorembodyfont{\rmfamily} \newtheorem{remark}{Remark}}
\newcommand{\abs}[1]{|#1|}
\title{Isothermic and S-Willmore Surfaces \\
as Solutions to Blaschke's Problem}
\author{Xiang Ma}
\begin{document}
\maketitle

\section*{\begin{center}Introduction\end{center}}
Isothermic surfaces and Willmore surfaces are important objects in 
differential geometry. Especially, they are surface classes invariant 
under M{\"o}bius transforms. Although seemed so distinct to each other, 
they may be introduced as the only non-trivial solutions to a problem
in the category of M{\"o}bius differential geometry(\cite{jeromin}).
\\

{\bf Blaschke's Problem:}
\textit{Let S be a sphere congruence with two envelops 
$f,\hat{f}:M^2\to S^3$, such that these envelops induce the same conformal
structure. Characterize such sphere congruences and envelop surfaces.}
\\
\\
Blaschke asked this question and solved it in \cite{blaschke}, namely 

\begin{theorem}
The non-trivial solution to the Blaschke's problem is either a pair of
isothermic surfaces forming Darboux transform to each other together with
the Ribaucour sphere congruence in between, or a pair of dual Willmore 
surfaces with their common mean curvature spheres. (Here \emph{non-trivial} 
means the two envelops are not congruent up to M{\"o}bius transforms.)
\end{theorem}

Note that Blaschke and his school only considered surfaces in 3-space and
ignored the higher co-dimension case. In $S^n$ this problem is still meaningful.
It is known that the construction of Darboux pair of isothermic surfaces as 
well as dual Willmore surfaces has a generalization to $S^n$
(see \cite{burstall,jeromin} and \cite{ejiri}), 
and they still constitute solutions to the generalized Blaschke's problem.
Here we will show they are exactly the only nontrivial solutions as before.
In 3-space every Willmore surface is S-Willmore, so our theorem reduces
to Blaschke's result. The main theorem reads as follows:

\begin{theorem}
If two distinct immersed surfaces into $S^n$ envelop a 2-sphere congruence 
and are conformal at correspondent points, then besides trivial cases, locally
it must be among either of the following two classes:\\
\hspace*{3mm} 1. a Darboux pair of isothermic surfaces;\\
\hspace*{3mm} 2. an S-Willmore surface with its dual surface.\\
In the first case, the correspondence between the two surfaces is orientation
preserving. In the second case, the correspondence is anti-conformal.
\end{theorem}

Naturally one may consider a special case of Blaschke's problem, that the 2-sphere
congruence is exactly the mean curvature spheres of either envelop surface.
The classification result is

\begin{theorem}\label{cmc}
Let $f$ be an surface immersed in $S^n$, whose mean curvature sphere congruence
has $\hat{f}$ as the second envelop with the same conformal structure. Then
either $f$ is M{\"o}bius equivalent to a cmc-1 surface in hyperbolic 3-space,
or $f$ is S-Willmore.
\end{theorem}

This paper is orgnized as follows. In section 1, we will recall the basic surface 
theory in M{\"o}bius geometry first, then treat the geometry of a pair of conformal
immersions, which is another contribution of this paper. We reformulate Blaschke's 
problem and divide the discussion into two cases, characterized as either $\theta$ 
or $\rho$ vanishes. It will be shown in section 2 that
non-trivial solutions in either cases are isothermic or S-Willmore surfaces
separately. Section 3 is devoted to the proof of Theorem~\ref{cmc}.
Finally, the geometric interpretation of the vanishness of $\theta$ or $\rho$ 
will be discussed in the appendix.

\section{Surface theory in conformal differential geometry}

\subsection{Light-cone model and fundamental equations}

In this part we follow the treatment in \cite{bpp}. Let $\mathbb{R}^{n+1,1}$
denote the Minkowski space, with isotropic lines correspondent to points
and space-like $(n-k)$-dim subspaces correspondent to $k$-spheres in $S^n$.

For a conformal immersion $f:M\to S^n$ of Riemann surface $M$, a 
\emph{(local) lift} is just a map $F$ from $M$ into the light cone such that
the null line spanned by $F(p)$ is $f(p)$. There is a decomposition
$M\times\mathbb{R}^{n+1,1} = V \oplus V^{\bot}$, where $V$ is a rank $4$ subbundle 
of $M\times\mathbb{R}^{n+1,1}$ of signature $(3,1)$ defined via local lift $F$
and complex coordinate $z$:
$$ V = \mbox{span} \{ F,dF,F_{z\bar{z}} \}.$$
\begin{remark}\label{mean}
Note $V$ is independent to the choice of local lift $F$, thus well-defined. It is
also M{\"o}bius invariant. The 2-sphere congruence given by $V$ (or equivalently 
the rank $n-2$ space-like subbundle $V^{\bot}$), called \emph{mean curvature 
sphere} or \emph{central sphere congruence}, is very important. The name 
\emph{mean curvature sphere} comes from the property that it is tangent to 
the surface and has the same mean curvature as the surface at the tangent 
point when the ambient space is endowed with a metric of Euclidean space 
(or any space form). 
\end{remark}
Given a local coordinate $z$. Among various choice of local lifts there is a 
canonical one, which is into the forward light cone and determined by 
$\abs{dF}^2=\abs{dz}^2.$
Such a canonical lift will be denoted by $Y$. A M{\"o}bius invariant framing 
of $V\otimes\mathbb{C}$ is given by 
$$\{Y,Y_z,Y_{\bar{z}},N\},$$
where we choose $N\in \Gamma(V)$ with
$\langle N,N\rangle=0,~\langle N,Y\rangle=-1,~\langle N,Y_z\rangle=0.$ 
These framing vectors are orthogonal to each other except that
$
\langle Y_{z},Y_{\bar{z}} \rangle = \frac{1}{2},~\langle Y,N \rangle = -1.
$
Especially, $\langle Y_{z},Y_{z}\rangle=0$ is the condition of conformality.

A fundamental equation in such a M{\"o}bius invariant surface theory is 
an inhomogeneous Hill's equation,
$$Y_{zz} + \frac{s}{2} Y = \kappa,$$
which defines two M{\"o}bius invariants, $s$ as a complex valued function, and 
$\kappa$ a section of $V^{\bot}\otimes\mathbb C$.
According to the observation in \cite{bpp}, $s$ is interpreted as the 
\emph{Schwarzian derivative} of immersion $f$, and $\kappa$ may be identified
with the normal valued \emph{Hopf differential}. When the local coordinate changes 
from $z$ to $w$, these quantities change according to the following rules,
\begin{eqnarray}
\kappa'&=&\kappa~
\Big(\frac{\partial z}{\partial w}\Big)^\frac{3}{2}
\Big(\frac{\partial\bar{z}}{\partial\bar{w}}\Big)^{-\frac{1}{2}},\label{dif-k}\\
s'&=& s~\Big(\frac{\partial z}{\partial w}\Big)^2+S_w(z).\label{dif-s}
\end{eqnarray}
where $S_w(z)$ is the classical Schwarzian derivative of $z$ with respect to $w$.

Finally, let $\psi\in\Gamma(V^{\bot})$ be an arbitrary section of the normal bundle
$V^{\bot}$, $D$ be the normal connection, we write down the structure equations,
\begin{equation} \label{mov-eq}
\left\{
\begin{array}{llll}
Y_{zz} &=& -\frac{s}{2} Y + \kappa, \\
Y_{z\bar{z}} &=& -\langle\kappa,\bar{\kappa}\rangle Y + \frac{1}{2}N, \\
N_{z} &=& -2 \langle\kappa,\bar{\kappa}\rangle Y_z - s Y_{\bar{z}}
+ 2D_{\bar{z}}\kappa, \\
\psi_z &=& D_z\psi + 2\langle\psi,D_{\bar{z}}\kappa\rangle Y
- 2\langle\psi,\kappa\rangle Y_{\bar{z}}.
\end{array}
\right.
\end{equation}
The conformal Gauss, Codazzi and Ricci equations as integrable conditions 
are 
\begin{eqnarray}
& \frac{1}{2}s_{\bar{z}} = 3\langle D_z\bar{\kappa},\kappa\rangle
+ \langle\bar{\kappa},D_z\kappa\rangle, \label{gauss}\\
& \mbox{Im}( D_{\bar{z}}D_{\bar{z}}\kappa
+ \frac{1}{2}\bar{s}\kappa ) = 0, \label{codazzi}\\
& R_{\bar{z}z}^D\psi
:= D_{\bar{z}}D_z\psi - D_z D_{\bar{z}}\psi
= 2\langle\psi,\kappa\rangle\bar{\kappa}
- 2\langle\psi,\bar{\kappa}\rangle\kappa. \label{ricci}
\end{eqnarray}

\subsection{Pair of conformal immersions}\label{pair}

Transforms of certain surface classes is an important topic in surface theory.
Generally for a given surface we construct the second immersion by solving 
certain differential equations. Such constructions usually admit kinds of duality 
between the old surface and the new one. In view of that, it is also natural to 
follow another line, namely to characterize such a pair of surfaces by some 
geometrical conditions. These considerations lead to a study of pair of conformal 
immersions, which seemed to be a rather interesting theme after dealing with 
single surfaces.

We start with two arbitrary conformal immersions $f,\hat{f}:M\to S^n$. 
Given coordinate $z$, to make a choice of normalized lifts, set $Y$ to be 
the canonical lift of $f$, and $\hat{Y}$, the corresponding lift of $\hat{f}$, 
to satisfy 
\begin{equation}\label{normal}
\langle Y,\hat{Y}\rangle=-1.
\end{equation}
With the canonical frames $\{Y,Y_z,Y_{\bar{z}},N\}$ at hand, we may express 
$\hat{Y}$ explicitly as
\begin{equation} \label{yhat}
\hat{Y}=\lambda Y + \bar{\mu}Y_z + {\mu}Y_{\bar{z}} + N + \xi,
\end{equation}
where $\lambda$ and $\mu$ are real-valued and complex-valued functions separately, 
and the real $\xi\in\Gamma(V^{\bot}).$ $\hat{Y}$ being isotropic implies 
$\lambda=\frac{1}{2}({\abs{\mu}}^2+\langle\xi,\xi\rangle).$
Substitute this back into (\ref{yhat}), and differentiate with respect to $z$,
a straightforward calculation with the help of (\ref{mov-eq}) yields the 
fundamental equation
\begin{equation} \label{yhat-z}
\hat{Y}_z =\frac{\mu}{2}\hat{Y}
+ \theta(Y_{\bar{z}} + \frac{\bar{\mu}}{2}Y)
+ \rho(Y_z + \frac{\mu}{2}Y) + \langle\xi,\zeta\rangle Y + \zeta ,
\end{equation}
where
\begin{eqnarray}
\theta &=& \mu_z - \frac{1}{2}{\mu}^2 - s - 2\langle\xi,\kappa\rangle, 
\label{theta} \\
\rho &=& \bar{\mu}_z - 2\langle\kappa,\bar{\kappa}\rangle
+ \frac{1}{2}\langle\xi,\xi\rangle, \label{rho} \\
\zeta &=& D_z\xi - \frac{\mu}{2}\xi +
2(D_{\bar{z}}\kappa + \frac{\bar{\mu}}{2}\kappa)~\in~\Gamma(V^{\bot}\otimes 
\mathbb{C}).
\end{eqnarray}

\begin{remark}\label{remark-geom}
What's the geometric meaning of $\theta$ and $\rho$? It is easy to check that 
$\theta$ and $\rho$ define a $(2,0)$ form and a $(1,1)$ form separately. 
We may define them with more freedom of choice, i.e.
\begin{eqnarray} 
\frac{\theta}{2}~{\mbox{d}z}^2 
&=& \frac{\langle Y_z\land Y,\hat{Y}_z\land\hat{Y}\rangle}
{\parallel Y\land\hat{Y}\parallel^2}~{\mbox{d}z}^2, \label{theta-geom}\\
\frac{\rho}{2}~\mbox{d}z\mbox{d}\bar z 
&=& \frac{\langle Y_{\bar z}\land Y,\hat{Y}_z\land\hat{Y}\rangle}
{\parallel Y\land\hat{Y}\parallel^2}~\mbox{d}z\mbox{d}\bar z.\label{rho-geom}
\end{eqnarray}
for arbitrary local lifts $Y,\hat{Y}$ and coordinate $z$, and independent to 
such choices. Here the inner product between bivectors are as usual. 
Note that if we interchange between $Y$ and $\hat{Y}$, $\rho$ turns to be 
$\bar\rho$, and $\theta$ invariant. They are invariants associated with a 
pair of immersed surfaces describing their relative positions. 
More discussion is left to the appendix.
\end{remark}

Now let us come down to Blaschke's problem. Assume $f$ and $\hat{f}$ form 
a pair of solution surfaces, $Y$ and $\hat{Y}$ the normalized lifts 
according to (\ref{normal}), $z=u+{\rm i}v$ the complex coordinate. 
The 2-sphere congruence tangent to $Y$ and passing through $\hat{Y}$ 
is given as span$\{ Y,{\rm d}Y,\hat{Y}\}$. According to the assumptions 
of Blaschke's problem, such 2-pheres also tangent to $\hat{Y}$, thus 
$\hat{Y}_{z}$ as well as $\zeta$ is contained in 
span$\{ Y,{\rm d}Y,\hat{Y}\}$. Since 
$\mbox{span}\{ Y,{\rm d}Y,\hat{Y}\}\oplus V^{\bot}$
is also a direct sum decomposition of $M\times\mathbb{R}^{n+1,1}$, 
this simply means
\begin{equation} \label{zeta}
\zeta=D_z\xi - \frac{\mu}{2}\xi + 2(D_{\bar{z}}\kappa + 
\frac{\bar{\mu}}{2}\kappa)=0.
\end{equation}
Under this circumstance, we have a simplified version of (\ref{yhat-z}):
\begin{equation} \label{yhat-z-new}
\hat{Y}_z =\frac{\mu}{2}\hat{Y}
+ \theta(Y_{\bar{z}} + \frac{\bar{\mu}}{2}Y)
+ \rho(Y_z + \frac{\mu}{2}Y).
\end{equation}
We calculate further that
$$
\langle\hat{Y}_z,\hat{Y}_z\rangle
=\langle\hat{Y}_z-\frac{\mu}{2}\hat{Y},\hat{Y}_z-\frac{\mu}{2}\hat{Y}\rangle 
=\theta\cdot\rho.
$$
By the last assumption that $\hat{Y}$ conformal to $Y$, we must have 
$\theta\cdot\rho=0$, i.e. $\theta=0$ or $\rho=0$. Thus Blaschke's problem 
is equivalent to finding all pairs of conformal immersions for which 
$\zeta=0$ and $\theta\cdot\rho=0$.

\section{Characterization of solutions}

\subsection{Preparation}

In this part, we will characterize any pair of surfaces solving Blaschke's 
problem. Let's choose the lifts $Y,\hat{Y}$ normalized as before, and assume 
$\zeta=0,~\theta\cdot\rho=0$. Note the first equation (\ref{zeta}) implies
\begin{eqnarray}
&&D_{\bar{z}}D_{\bar{z}}\kappa + \frac{1}{2}\bar{s}\kappa \label{kappa}\\
&=& D_{\bar{z}}(-\frac{1}{2}D_z\xi + \frac{\mu}{4}\xi - 
\frac{\bar{\mu}}{2}\kappa)+ \frac{\bar{s}}{2}\kappa \nonumber\\
&=& -\frac{1}{2}D_{\bar{z}}D_z\xi + \frac{\mu}{4}D_{\bar{z}}\xi + 
\frac{\mu_{\bar{z}}}{4}\xi
- \frac{\bar{\mu}}{2}D_{\bar{z}}\kappa - \frac{\bar{\mu}_{\bar{z}}}{2}\kappa
+ \frac{\bar{s}}{2}\kappa \nonumber\\
&=& -\frac{1}{2}D_{\bar{z}}D_z\xi + \frac{\mu}{4}D_{\bar{z}}\xi + 
\frac{\mu_{\bar{z}}}{4}\xi
- \frac{\bar{\mu}}{2}(-\frac{1}{2}D_z\xi + \frac{\mu}{4}\xi - 
\frac{\bar{\mu}}{2}\kappa)
- \frac{\bar{\mu}_{\bar{z}}}{2}\kappa + \frac{\bar{s}}{2}\kappa \nonumber\\
&=& \Big(\!-\frac{1}{2}D_{\bar{z}}D_z\xi-\langle\xi,\bar{\kappa}\rangle\kappa\Big)
+ \Big(\frac{\mu}{4}D_{\bar{z}}\xi + \frac{\bar{\mu}}{4}D_{z}\xi \Big)
- \frac{1}{8}{\abs{\mu}}^2\xi
+ \frac{\mu_{\bar{z}}}{4}\xi - \frac{\bar{\theta}}{2}\kappa. \nonumber
\end{eqnarray}
Codazzi and Ricci equations (\ref{codazzi})(\ref{ricci}) now implies
\begin{equation}\label{real}
\mbox{Im}(\mu_{\bar z}\xi-2\bar\theta\kappa)=0
\end{equation}
On the other hand, $\zeta=0$ together with Gauss equation (\ref{gauss}) yields
\begin{eqnarray}
\theta_{\bar{z}}
&=& (\mu_{\bar{z}})_z - \mu {\mu_{\bar{z}}} - s_{\bar{z}}
- 2{\langle\xi,\kappa\rangle}_{\bar{z}} \label{theta-zbar}\\
&=&
\Big(\bar\rho+2\langle\kappa,\bar{\kappa}\rangle
-\frac{1}{2}\langle\xi,\xi\rangle\Big)_z 
-\mu\Big(\bar\rho+2\langle\kappa,\bar{\kappa}\rangle
-\frac{1}{2}\langle\xi,\xi\rangle\Big) \nonumber\\
&& -6\langle D_z\bar{\kappa},\kappa\rangle-2\langle\bar{\kappa},D_z\kappa\rangle 
-2{\langle\xi,\kappa\rangle}_{\bar{z}} \nonumber\\
&=& \bar\rho_z-\mu\bar\rho 
-2\langle D_{\bar{z}}\xi+2D_z{\bar\kappa}+\mu\bar{\kappa},\kappa\rangle
+\langle\xi,-D_z\xi+\frac{\mu}{2}\xi-2D_{\bar{z}}\kappa\rangle \nonumber\\
&=& \bar\rho_z-\mu\bar\rho-\langle\bar{\mu}\xi,\kappa\rangle 
+\langle\xi,\bar{\mu}\kappa\rangle \nonumber\\
&=& \bar\rho_z-\mu\bar\rho. \nonumber
\end{eqnarray}
These formulae will be very useful in our discussion. Also note that if
$\theta=\rho=0$ on an open subset of $M$, then $\hat{Y}_z=\frac{\mu}{2}\hat{Y}$,
implying that $\hat{Y}$ correspond to a fixed point in $S^n$. 
After a stereographic projection from this point, it turns out to be 
an immersion into $\mathbb{R}^n$ with its tangent planes passing through 
$\infty$. Such a trivial and degenerate case is excluded from our consideration.

\subsection{{\boldmath $\rho = 0\neq\theta$}: Isothermic case}

Let's consider the first case, $\rho = 0\neq\theta$.
Note that $\rho=0$ implies $\bar{\mu}_z$ is real-valued by (\ref{rho}).
By (\ref{real}), $\bar\theta\kappa$ has to be real, too. With $\rho=0$, 
(\ref{theta-zbar}) implies $\theta$ is holomorphic.
So $\kappa=\theta\cdot\frac{\bar\theta\kappa}{\abs{\theta}^2}$, where 
the first term is holomorphic, and the next term is real. We may define 
a new holomorphic coordinate $w$ by ${\rm d}w = \sqrt{\theta}~{\rm d}z$. 
Due to the transformation rule of $\kappa$ given by (\ref{dif-k}), we get 
${\kappa'}=\frac{\bar\theta\kappa}
{\displaystyle \abs{\theta}^{3/2}}$, a real vector-valued form.

Recall that classically, a surface in $\mathbb{R}^3$ is called 
\emph{isothermic} if it can be conformally parametrized by its curvature 
lines. This notion is indeed conformal invariant. It has been generalized 
(\cite{burstall,palmer,schief}) to higher codimension spaces by an 
equivalent characterization that the Hopf differential is real valued 
in a suitable complex coordinate. So our analysis shows $f$ must be 
isothermic in this case. 

Furthermore, if we consider $\hat{f}$ as the first (original) one in such 
a pair, on account of (\ref{theta-geom}) and (\ref{rho-geom}), we have the 
same invariant $\rho$ up to complex conjugation, and the same $\theta$. 
Similarly $\bar\theta\hat\kappa$ is real, and when the coordinate $z$ is 
chosen such that $\kappa$ is real, $\hat\kappa$ will be real at the same 
time. That means both surfaces are isothermic and their curvature lines 
correspond. By the characterization given in \cite{burstall,jeromin}, 
such two surfaces are said to envelop a conformal Ribaubour sphere 
congruence, and they form Darboux transform to each other. Especially, 
given such an coordinate $z$ and real $\kappa$, the holomorphic $\theta$ 
must be also real-valued, thus constant. This constant may be identified 
with the real parameter appearing in the construction of Darboux transforms 
up to a choice of certain holomorphic 2-form over the underlying $M$ 
(compared with \cite{burstall,jeromin-musso}). 

Note that after close examination, the definition of Darboux transform of an 
isothermic surface given in \cite[${\S}$ 2.2.2]{burstall} may be regard
as a purely geometric characterization that \textit{such a pair envelop a 
2-sphere congruence with conformal metric and the same orientation induced 
by those 2-spheres}. (As to the orientation issue, we will clarify in the 
appendix.) Regretfully, both Burstall and Hertrich-Jeromin only mentioned the 
characterization assuming the sphere congruence to be 
\emph{Ribaucour} and overlooked our simpler version.

\subsection{{\boldmath $\theta= 0\neq\rho$}: S-Willmore case}

The second possibility is $\theta = 0 = \rho$ on an open subset.
In this part, we first consider the seemingly \emph{trivial} case, $\xi= 0$.
Put this into (\ref{kappa}), we see immediately that
\begin{equation}\label{willmore}
D_{\bar{z}}D_{\bar{z}}\kappa + \frac{1}{2}\bar{s}\kappa = 0.
\end{equation}
In \cite{bpp}, the authors proved that equation (\ref{willmore}) characterizes
\emph{Willmore} surfaces. Interest in such surfaces originated long ago. 
The reader may consult \cite[Ch.7]{willmore} and \cite[Ch.3]{jeromin}
for a historical retrospect and references therein. We should emphasize here 
that $Y$ satisfies a condition stronger than Willmore, namely, $\xi=0$ and 
$\zeta=0$ implies
\begin{equation}\label{strong}
D_{\bar{z}}\kappa + \frac{\bar{\mu}}{2}\kappa = 0.
\end{equation}
\begin{definition}
A Willmore surface is called \emph{S-Willmore} if its Hopf differential 
$\kappa$ is linear dependent to $D_{\bar{z}}\kappa$.
\end{definition}
\begin{remark}
Note the property of being S-Willmore is well-defined, i.e. independent to 
the choice of $z$. It may be easily checked by formula~(\ref{dif-k}).
\end{remark}

The concept of a \emph{S-Willmore} surface has appeared in \cite{wang} as 
\emph{strong Willmore} surface. More early, it was given by Ejiri 
(\cite{ejiri}) with some slight modification. Namely, since (\ref{strong}) 
holds automatically for surfaces in 3-space, any Willmore surface in $S^3$ 
is S-Willmore, which is excluded from Ejiri's definition. 
In our opinion, including the codim-1 case into the definition is more 
natural and favorable, because they share the same \emph{duality} and 
arise as one important class of solutions to Blaschke's problem. Such a 
duality has attracted many geometers (see \cite{bryant,ejiri,wang}). 
Let's briefly review it at here.

Suppose $Y$ is the canonical lift of a S-Willmore surface. Locally we may 
assume it has no umbilic points, i.e. $\kappa\neq 0$. Thus (\ref{strong})
holds for some $\mu$. We construct a second surface via this $\mu$ and
$\hat{Y}=\frac{\abs\mu^2}{2}Y+\bar{\mu}Y_z+{\mu}Y_{\bar{z}} + N,$
i.e. equation (\ref{yhat}) with $\xi=0$. At this moment we see the invariant
$\theta=0$ directly due to the Willmore condition (\ref{willmore}) and 
$\xi=0$. Our choice of $\mu$ and $\xi$ also yields $\zeta=0$. 
So (\ref{yhat-z}) is reduced to  
$\hat{Y}_z =\frac{\mu}{2}\hat{Y}+\rho(Y_z + \frac{\mu}{2}Y).$
This implies 
$\mbox{span}\{ Y,Y_{z},Y_{\bar{z}},Y_{z\bar{z}}\}
= \mbox{span}\{ \hat{Y},\hat{Y}_{z},\hat{Y}_{\bar{z}},\hat{Y}_{z\bar{z}}\}$,
so $Y$ and $\hat{Y}$ share the same mean curvature sphere (yet with opposite
orientation, see appendix). Such a duality enables us to conclude that 
$\hat{Y}$ is also S-Willmore.

\begin{remark}
General Willmore surfaces may not be S-Willmore, and we can hardly expect 
any duality results for them. For such examples see the remark in the final 
section of \cite{ejiri}
and references therein.
\end{remark}

\subsection{Trivial case}\label{trivial}

Finally we come to the discussion of the {\sl non-trivial} part of case
$\theta= 0\neq\rho$, the occasion that $\xi\neq 0$. This time (\ref{real}) 
implies that $\bar\mu_z$ as well as $\rho$ must be real. Together with 
(\ref{theta-zbar}) we obtain
$\rho_z=\mu\rho.$
Since $\rho\neq 0$, without loss of generality we may assume it to be positive
and use it to scale $\hat{Y}$. Let $\tilde{Y}=\frac{1}{\rho}\hat{Y}$ 
(This is indeed the canonical lift of $\hat{f}$). 
Define $X:=\tilde{Y}-Y$ which is real-valued.
Then $\rho_z=\mu\rho, \theta=0$ and (\ref{yhat-z-new}) yields
\begin{eqnarray*}
X_z=(\frac{1}{\rho}\hat{Y}-Y)_z
=-\frac{\mu}{\rho}\hat{Y} + \frac{1}{\rho}\Big(\frac{\mu}{2}\hat{Y}
+\rho(Y_z + \frac{\mu}{2}Y)\Big)-Y_z
=-\frac{\mu}{2}X,
\end{eqnarray*}
Thus the real line spanned by $X$ is constant. We normalize $X$ by defining
$\tilde{X} = \sqrt{\rho}\cdot X$. Then we will have $\tilde{X}_z = 0$.
Thus $\tilde{X}$ is a fixed vector in $R^{n+1,1}$. Observe that
$$
\begin{array}{l l}
\langle X,X \rangle
= \langle \frac{1}{\displaystyle \rho}\hat{Y}-Y,\frac{1}{\displaystyle 
\rho}\hat{Y}-Y \rangle
= \frac{2}{\displaystyle \rho}, \\[1mm]
\langle Y,X \rangle
= \langle Y,\frac{1}{\displaystyle \rho}\hat{Y}-Y \rangle
= -\frac{1}{\displaystyle \rho},
\end{array}
$$
we find
$$
\tilde{Y}=Y+X=Y-\frac{2\langle Y,X \rangle}{\langle X,X \rangle}X
=Y-\frac{2\langle Y,\tilde{X}\rangle}{\langle\tilde{X},\tilde{X}\rangle}
\tilde{X},
$$
which means that $\tilde{Y}$ is a reflection of $Y$ with respect to 
$\tilde{X}$. So the underlying map $\hat{f}$ differs with $f$ only by 
a M\"{o}bius transform. Such a trivial case is easy to understand.

\subsection{Conclusion and comments}

Sum together, our main theorem is proved. 

We are interested in such a characterization not only because this problem is 
natural and interesting in itself, but also due to our concern about transforms
of surfaces. Note there are already many beautiful results about transformations 
of isothermic surfaces, e.g. dual isothermic surfaces, Bianchi and Darboux 
transforms together with the permutability theorems (\cite{burstall,jeromin}). 
As to Willmore surfaces, there are only some partial results 
(\cite{quater}), which seemed similar to isothermic case yet more subtle.
We think the relationship between these two surface classes is worthy 
of more exploration. 
 
In the appendix, $\rho=0$ is interpreted as \emph{touching} condition, 
whereas $\theta=0$ is \emph{co-touching}. At here they separately indicate 
that in the isothermic case the correspondence is orientation preserving 
and the Willmore case orientation reversing. In view that $\rho=0$ characterize
the Darboux transforms of isothermic surfaces, one would suspect that the
condition $\theta=0$ might be characteristic of transforms of Willmore 
surfaces. In other words, for a given Willmore surface in $S^n$, one may 
consider a second conformally immersion with $\theta=0$, which turns out
to be also Willmore under additional assumptions. This is the key idea in 
the construction of adjoint transforms of Willmore surfaces (\cite{ma}).
(In that work the mean curvature sphere of the original surface is assumed
to pass through the second one, yet not neccessarily tangent to it.)
Why such an idea works is still mysterious to the author. Here we only
remark that Blaschke's problem provokes further thought on these classical 
objects and shed new light on the construction of transforms of Willmore 
surfaces.

\section{Envelopped mean curvature spheres}

In the category of M{\"o}bius geometry, we may further ask: when does such 
a 2-sphere congruence solving Blaschke's problem happen to be the
mean curvature sphere of either envelop surface? The answer is
Theorem~\ref{cmc}, and here we will give the detail.

{\flushleft\it Proof to Theorem~\ref{cmc}:}\quad
We have seen that one class of previous solutions, dual S-Willmore 
surfaces, envelopping the same mean curvature spheres, also satisfy 
the further assumption. 

Now consider the isothermic case. 
Suppose $f$ is an isothermic surface with Darboux transform $\hat{f}$, 
both envelopping the mean curvature spheres of $f$. Then the lift $\hat Y$ 
as well as $\hat{Y}_z$ must be contained in 
$V=\mbox{span}\{Y,Y_z,Y_{\bar z},N\}$. By (\ref{yhat})(\ref{yhat-z}),
$\xi=0=\zeta$, hence $D_{\bar z}\kappa + \frac{\bar\mu}{2}\kappa=0$,
i.e. $f$ satisfies (\ref{strong}). Indeed, (\ref{strong}) holds if 
and only if the mean curvature spheres of $f$ has a second envelop
(\cite{wang}). Without loss of generality we may assume that locally
$\kappa$ and $\theta$ is non-zero, and choose the complex coordinate 
$z=u+{\rm i}v$ so that $\kappa$ is real. we immediately find that the 
real frames $\{Y,Y_u,Y_v,N,\kappa\}$ satisfy a PDE system given by 
(\ref{mov-eq}) and (\ref{strong}). Integration shows that $Y$ is 
contained 
in a 5-dim Minkowski subspace. So $f$ is an immersion into $S^3$. 
As the second envelop surface of the mean curvature spheres of $f$, 
$\hat{f}$ must also be contained in the same 3-space. The canonical 
lift of $\hat{f}$ is given via $\tilde{Y}=\frac{1}{\theta}\hat{Y}$ 
due to (\ref{yhat-z-new}) and $\rho=0$. Keeping in mind that $\theta$
is real and holomorphic, hence constant, and that 
$\bar{\mu}_z - 2\langle\kappa,\bar{\kappa}\rangle
=\rho-\frac{1}{2}\langle\xi,\xi\rangle=0$, we differentiate $\tilde{Y}$:
\begin{eqnarray*}
\tilde{Y}_z&=&\frac{1}{\theta}\hat{Y}_z=\frac{\mu}{2}\tilde{Y}
+(Y_{\bar z}+\frac{\bar\mu}{2}Y),\\
\tilde{Y}_{zz}&=&\frac{\mu}{2}\tilde{Y}_z+\frac{\mu_z}{2}\tilde{Y}
+(Y_{\bar{z}z}+\frac{\bar\mu}{2}Y_z+\frac{\bar\mu_z}{2}Y)\\
&=&\frac{\mu}{2}\Big[\frac{\mu}{2}\tilde{Y}
+(Y_{\bar z}+\frac{\bar\mu}{2}Y)\Big]+\frac{\mu_z}{2}\tilde{Y}
+(\frac{1}{2}\bar{\mu}_z-\langle\kappa,\bar{\kappa}\rangle)Y
+\frac{\bar\mu}{2}Y_z+\frac{1}{2}N\\
&=&(\frac{\mu_z}{2}+\frac{\mu^2}{4})\tilde{Y}+\frac{1}{2}
(\frac{1}{2}\abs{\mu}^2 Y +\bar\mu Y_z + \mu Y_{\bar z} + N)\\
&=&(\cdots)\tilde{Y}.
\end{eqnarray*}
By definition $\hat\kappa=0$, and $\hat{f}$ is a round sphere. 
If we restrict all these objects in 3-space, and endow
the 3-ball enclosed by $\hat{f}$ with the hyperbolic metric of constant
curvature -1, then $\hat{f}$ stands as the boundary at infinity of this 
hyperbolic space. Each mean curvature sphere of $f$ is tangent 
to this boundary 2-sphere, thus a horo-sphere with curvature 1. 
As pointed out in Remark~\ref{mean}, that implies $f$ is of constant
mean curvature 1.

Note in the degenerate case that $\theta=\rho=0$, namely a surface in 
$\mathbb{R}^n$ with all its tangent planes (passing through $\infty$),
if these tangent planes are mean curvature spheres of this surface,
then it is a minimal surface in this Euclidean space, which is regarded
as a degenerate case of S-Willmore. The trivial case in ${\S}$
\ref{trivial} is readily excluded since $\xi\neq 0$, hence the mean 
curvature spheres of $f$ never pass through $\hat{f}$.
\quad$\Box$

\begin{remark}
In fact we have obtained a characterization of cmc-1 surfaces in 
hyperbolic 3-space as the only surfaces in $S^n$ whose mean curvature
spheres has a second envelop surface and they correspond conformally 
with compatible orientation. Such a result first appeared in 
\cite{jeromin-musso}, which characterizes these surfaces among 
isothermic surfaces in 3-space as those ones whose mean curvature 
spheres give rise to a second envelop surface and they form a 
Darboux pair. Our generalization is in accord with the generalization 
of Blaschke's problem. Interestingly that general cmc surfaces in 
hyperbolic $n$-space do not share the same property, mainly due to 
the fact that in higher codimension case they fail to be isothermic.
\end{remark}

\section*{Appendix: touching and co-touching}

As a sub-geometry of Lie-sphere geometry, M{\"o}bius geometry validates
the notion of \emph{contact elements}. Here we always 
assume them to be oriented. Such a 2-dim contact element is just a 
2-dim oriented subspace in the tangent space of $S^n$ at one point $p$, 
which corresponds to a 3-dim oriented subspace of signature $(2,0)$
in the Minkowski space $\mathbb{R}^{n+1,1}$. In the discussion below
we will always represent such a contact element as frames $\{X,X_1,X_2\}$, 
where $X$ is light-like and span the null line given by $p$, $\{X_1,X_2\}$ 
an oriented conformal (i.e. $\langle X_1,X_2\rangle=0,
\langle X_1,X_1\rangle=\langle X_2,X_2\rangle>0$) frame in the orthogonal 
complement of $X$. We also utilize the \emph{complex contact element}
represented by $X\land(X_1-{\rm i}X_2)$. The latter is also represented
as $X\land X_z$ where $X_z=\frac{\partial}{\partial z}X$ and 
$z=u+{\rm i}v$ is a complex coordinate.

Consider two contact elements $\Pi=\{Y,Y_1,Y_2\}$ and $\hat\Pi=\{\hat Y,
\hat Y_1,\hat Y_2\}$. The null lines spanned separately by $Y,\hat Y$ 
are distinct, so these contact elements are at different points $p,\hat p$. 
Inspired by the formulae (\ref{theta-geom})(\ref{rho-geom}) in 
Remark~\ref{remark-geom}, we define two invariants associated with such 
a pair of 2-dim contact elements:
\begin{eqnarray} 
\frac{\theta}{2}&=& \frac{\langle Y\land(Y_1-{\rm i}Y_2),
\hat{Y}\land(\hat{Y}_1-{\rm i}\hat{Y}_2)\rangle}
{\parallel Y\land\hat{Y}\parallel^2}, \label{theta-contact}\\
\frac{\rho}{2}&=& \frac{\langle Y\land(Y_1+{\rm i}Y_2),
\hat{Y}\land(\hat{Y}_1-{\rm i}\hat{Y}_2)\rangle}
{\parallel Y\land\hat{Y}\parallel^2}.\label{rho-contact}
\end{eqnarray}
Note they are determined up to the choice of conformal frames(coordinates)
at correpondent points. The M{\"o}bius invariance is obvious.

How about two contact elements at the same point $p$? Suppose they are 
given as $\{Y,Y_1,Y_2\}$ and $\{Y,\hat Y_1,\hat Y_2\}$. We need only to
consider the 2-planes $\mbox{span}\{Y_1,Y_2\}$ and 
$\mbox{span}\{\hat Y_1,\hat Y_2\}$. The inner products between these
frame vectors yields a $2\times 2$ matrix 
$A=(a_{ij})=\Big(\langle Y_i,\hat Y_j\rangle\Big)$
with decomposition
$$
A=Q+P=\frac{1}{2}
\left(\begin{array}{cc}
a_{11}+a_{22} & a_{12}-a_{21} \\
a_{21}-a_{12} & a_{11}+a_{22} 
\end{array}\right)
+\frac{1}{2}\left(\begin{array}{cc}
a_{11}-a_{22} & a_{12}+a_{21} \\
a_{12}+a_{21} & a_{22}-a_{11}
\end{array}\right).
$$
This may be regarded as a decomposition into the commutative and 
anti-commutative parts with respect to 
$$
J=\left(\begin{array}{rc}
0 & 1 \\ -1 & 0
\end{array}\right)
$$
i.e. $Q=\frac{1}{2}(A-JAJ),P=\frac{1}{2}(A-JAJ)$.
Note that $A$ is determined up to multiplication from left and
right by orthogonal matrices and (nonzero) scalar matrices, 
and such a decomposition is invariant up to these actions. 
Corresponding to $Q$ and $P$, there are two complex numbers
$$\begin{array}{ccccc}
\theta'&=&(a_{11}+a_{22})+{\rm i}(a_{12}-a_{21})
&=&\langle Y_1+{\rm i}Y_2,\hat Y_1-{\rm i}\hat Y_2\rangle,\\
\rho'&=&(a_{11}-a_{22})+{\rm i}(a_{12}+a_{21})
&=&\langle Y_1-{\rm i}Y_2,\hat Y_1-{\rm i}\hat Y_2\rangle,
\end{array}$$
which are analogues as previously defined $\theta$ and $\rho$.
Now we can introduce

\begin{definition}
Two contact elements $\Pi_1$ and $\Pi_2$ at one point are said to
\emph{touch} each other if $\rho'=0$ and \emph{co-touch} each 
other if $\theta'=0$.

Two oriented conformally immersed surfaces intersecting at $p$
\emph{touch~(co-touch)} each other if the contact elements given
by their tangent spaces at $p$ touch~(co-touch).
\end{definition}

\begin{remark}
As the author knows, the notion of {\it left/right touch} was first 
considered by Pedit and Pinkall and may be found in \cite{bohler}.
Restricted to the case of 4-space, we may identify $\mathbb{R}^4$ with
quaternions $\mathbb{H}$ and define \emph{left and right normal vectors} 
for any oriented 2-plane, which is an algebraic version of the well-known 
fact that $Gr_+(2,4)\simeq S^2\times S^2.$ If the tangent planes of two 
surfaces at the intersection point share the same left (right) normal 
vector, then they are said to \emph{touch each other from left (right)}.
If the touching is both from left and right, then they are tangent with
the same induced orientation. In the quaternionic framework they found
this M{\"o}bius invariant notion and used it to define \emph{Darboux 
transforms of general surfaces} which is a generalizaton of the classical
Darboux transforms of isothermic surfaces (\cite{bohler}). Following their
observation, in \cite{quater} it is shown that the mean curvature spheres
of a Willmore surface touches its \emph{forward (backward) two-step 
B{\"a}cklund transforms} from left (right) yet with a negative sign. This
motived the author to introduce the notion of \emph{co-touch}. Here we 
introduce them in a slightly different yet equivalent way.
\end{remark}

\begin{remark}
The duality between \textit{left} and \textit{right} is only due to the 
orientation induced by the identification $\mathbb{R}^4=\mathbb{H}$, 
hence not essencial. Put this aside, the notion of touching and 
co-touching are independent to the choice of such an identification.
On the other hand, if either surface reverse its orientation, touching
will change to be co-touching, and vice versa. In case that the 2-planes
are in the same 2-space, it is easy to verify that touching corresponds to
compatible orientations and co-touching corresponds to opposite ones.
\end{remark}

\begin{remark}
In this work we avoid the use of quaternions and shows the key point is
the product matrix $A$ defined for any pair of 2-planes. Such matrices
are considered to be equivalent up to multiplication by orthogonal and 
nonzero scalar matrices, and the classification under this equivalence
relationship naturally leads to our invariants. For pairs of $m$-dim
subspaces also exist similar conformal invariants. Yet the best case is
still surfaces due to the underlying complex structure.
\end{remark}

To clarify the geometric meaning of $\theta=0$ and $\rho=0$, observe 
that given coordinate $z=u+{\rm i}v$, contact element $\Pi_1=\{Y,Y_u,Y_v\}$ 
at $Y(p)$, and single point $\hat Y(p)$, there is an unique oriented 
2-sphere passing through $Y(p),\hat Y(p)$ and tangent to $Y$ 
with compatible orientation. It is given by the 4-dim subspace of 
signature $(3,1)$ spanned by $\{Y,Y_u,Y_v,\hat Y\}$, with the orientation 
fixed by the oriented contact element $\Pi=\{Y,Y_u,Y_v\}$ or the 
complexification $Y\land Y_z$. Denote it as $S(p)$. Now we may state

\begin{proposition}
Given two conformal immersions $f,\hat f$, the invariant $\rho(p)=0$ iff 
the 2-sphere $S(p)$ touches $\hat Y$ at $p$, and $\theta(p)=0$ iff 
$S(p)$ co-touches $\hat Y$ at $p$.
\end{proposition}

{\flushleft\it Proof}\quad
Since we have the freedom of choice of lifts, we may take the normalized
lifts $Y,\hat Y$ as in subsection~\ref{pair}. Here 
$\hat{Y}=\frac{\abs\mu^2}{2}Y+\bar{\mu}Y_z+{\mu}Y_{\bar{z}} + N$ 
is orthogonal to $Y_z + \frac{\mu}{2}Y$. Note that under the reflection 
with respect to $Y-\hat Y$, $S(p)$ is invariant with reversed orientation, 
and the complex contact element $\Pi_1=Y\land(Y_z+\frac{\mu}{2}Y)$
is mapped to $\hat Y\land(Y_z+\frac{\mu}{2}Y)$. Thus the complex contact 
element given by $S(p)=\mbox{span}\{Y,Y_u,Y_v,\hat Y\}$ at $\hat Y(p)$ 
should be $\Pi'=\hat Y\land(Y_{\bar z}+\frac{\bar\mu}{2}Y)$. 
On the other hand, the complex contact element given by immersion 
$\hat Y$ at $p$ is $\Pi_2=\hat Y\land \hat Y_z$. 
By definition and the fundamental equation (\ref{yhat-z}), we see 
$$
\theta'=\langle Y_z+\frac{\mu}{2}Y,\hat Y_z\rangle=\theta,~~
\rho'=\langle Y_{\bar z}+\frac{\bar\mu}{2}Y,\hat Y_z\rangle=\rho.
$$
The conclusion now follows from the definition of touching and co-touching.
\quad$\Box$

\begin{corollary}
As solution to Blaschke problem, the surface pair share compatible 
orientations at correspondent points in the isothermic case ($\rho=0$),
thus a Darbux pair, or inducing opposite orientations on the common
mean curvature spheres in the S-Willmore case ($\theta=0$). 
\end{corollary}

\end{document}